\documentclass{amsart}
\usepackage{graphicx}
\usepackage{amssymb, amsthm, amsmath}
\usepackage{tikz}
\usetikzlibrary{arrows}
\usetikzlibrary{decorations.markings}
\newtheorem{theorem}{Theorem}[section]
\newtheorem*{theorem*}{Theorem}

\newtheorem{lemma}[theorem]{Lemma}
\newtheorem{prop}[theorem]{Proposition}
\newtheorem{fact}[theorem]{Fact}
\newtheorem{claim}[theorem]{Claim}
\newtheorem{example}[theorem]{Example}

\theoremstyle{definition}
\newtheorem{defn}[theorem]{Definition}
\theoremstyle{remark}
\newtheorem{rem}[theorem]{Remark}
\numberwithin{equation}{section}
\usepackage{enumerate}

\newcommand{\sub}{\subseteq}
\newcommand{\CM}{{\mathcal M}}
\newcommand{\CL}{{\mathcal L}}

\newcommand{\la}{\langle}
\newcommand{\ra}{\rangle}
\newcommand{\CN}{\mathcal N}
\newcommand{\Loc}[2]{\mu({#1}/{#2})}

\usepackage{lipsum}
\usepackage{array}
\title{An o-minimalist view of the group configuration}
\author{Ya'acov Peterzil}

\address{kobi@math.haifa.ac.il}
\email{}%

\thanks{This note was written for the book ``Research Trends in Contemporary Logic'' of College Publications}

\begin{document}

\begin{abstract}
The group configuration in o-minimal structures gives rise, just like in the stable case, to a transitive action of  a type-definable group  on a partial type. Because $acl=dcl$  the o-minimal proof is significantly simpler than Hrushovski's original argument. Several equivalent versions, which are more suitable to the o-minimal setting, are formulated, in functional language and also in terms of a certain $4$-ary relation.

In addition, the following question is considered: Can every definably connected type-definable group be definably embedded into a definable group of the same dimension? Two simple cases with a positive answer are given.

\end{abstract}

\maketitle

\section{Introduction}

The group configuration is one of the most important tools of Geometric Model Theory.
It allows to extract a (type) definable group from a small number of  model theoretic (in)dependencies.
Logicians working in o-minimal structures are often asked by other model theorists whether the group configuration theorem holds in o-minimal structures.  My usual reply is: ``Yes, but it does not suit well the o-minimal setting''.
In fact, even in \cite{trichotomy}, where we recovered a group and a field in the o-minimal setting, we did not make use
of such a configuration.
In this note I will try to clarify this point, proving a precise result about the recovery of a (type definable) group and its action from an $acl$-group configuration and then discussing some equivalences.

Roughly described, the proof of Hrushovski's Group Configuration theorem goes via three main steps. In the first  one replaces the original configuration where the dependencies are determined by {\em algebraic closure}, by an equivalent configuration where the dependencies  are given by {\em definable closure}. In the second step  one uses ``germs'' of definable functions to produce  a ``group chunk'' from this $dcl$-configuration and in the final step a (type) definable group is obtained using the so-called Hrushovski's group chunk theorem, a generalization of a theorem of Weil from the theory of algebraic groups.

The proof of the o-minimal result skips the above first step, since in linearly ordered structures $acl=dcl$. This simplifies substantially the argument. Then, using the notion of {\em an infinitesimal locus} and definable functions on them one obtains from the configuration an associativity-like functional equation.  As in the stable proof, composition of  functions yield a type-definable group and action, on infinitesimal types.

Before stating the main theorem recall that {\em a type-definable group} is a partial type $G(x)$ (i.e. a collection of formulas in $x$ of cardinality at most $\kappa=|L|$) together with a definable  binary function whose restriction to $G\times G$ is a group operation. A group action of $G$ on a type $\Sigma$ is called {\em definable} if it is the restriction to $G\times \Sigma$ of a definable function.

  \begin{theorem*} Let $\CM$ be an $|L|^+$-saturated o-minimal structure which eliminates imaginaries.
   Assume that $V=(a_1,a_2,a_3,x_1,x_2,x_3)$ is an $(m,k)$-configuration over $A$ (see \ref{defn1} for precise definition).

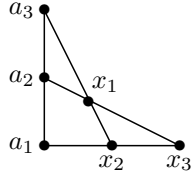
\begin{equation}
  \begin{tikzpicture}[scale=0.9]
\draw[black,line  width=0.2mm] (0,0) -- (0,2);
\draw[black,line  width=0.2mm] (0,0) -- (2,0);
\draw[black,line  width=0.2mm] (1,0) -- (0,2);
\draw[black,line  width=0.2mm] (0,1) -- (2,0);
\node at (0,0) {$\bullet$};
\node[left,scale=1.0] at (0,0) {$a_1$};
\node at (1,0) {$\bullet$};
\node[below,scale=1.0] at (1,0) {$x_2$};
\node at (2,0) {$\bullet$};
\node[below,scale=1.0] at (2,0) {$x_3$};
\node at (0,1) {$\bullet$};
\node[left,scale=1.0] at (0,1) {$a_2$};
\node at (0,2) {$\bullet$};
\node[left,scale=1.0] at (0,2) {$a_3$};
\node at (0.65,0.65) {$\bullet$};
\node at (0.9,0.9) {$x_1$};
\end{tikzpicture}
\end{equation}

 Then there exists an $m'$-dimensional type-definable group $(G,\star)$ over $M$, with $k\leq m'\leq m$,
acting definably and transitively on a $k$-dimensional partial type over $M$.

The configuration and group action are related as follows:  There is $B\supseteq A$ independent from $V$ over $A$, and there are $g_1,g_2,g_3\in G$, and $y_1,y_2,y_3 \models \Sigma$ such that for each $i=\{1,2,3\}$, we have $g_i\in dcl(a_iB)$,  and $dcl(x_iB)=dcl(y_iB)$, and in addition $(g_1,g_2,g_3,y_1,y_2,y_3)$ is an $(m',k)$-configuration in $ G$ as in Example \ref{gpexample}.

If, in addition, the configuration is minimal then $m'=m$ and each $g_i$ is inter-definable with $a_i$ over $B$.
\end{theorem*}

As in the stable case, an additional node in the configuration (see \cite[Appendix]{Bays-Hils-Moosa}) yields an abelian group. Also, if $m=k$ then one can identify the group with the homogeneous space.

\begin{rem} \begin{enumerate}
\item The assumption that $\CM$ eliminates imaginaries is considered mild in the o-minimal setting. E.g., it holds if $\CM$ expands an ordered group.
In fact, using the Trichotomy Theorem it can be eliminated but it is here in order to avoid technicalities.
\item The domain of the group $G$ is the infinitesimal type of $a_1$, or more precisely its quotient  by a definable equivalence relation. It acts on the infinitesimal type of $x_2$ or any other $x_i$. These are types over the model $\CM$, but after the fact could be chosen over any model which contains the configuration and the parameter set $B$ above.
\item In Section \ref{quadrangle} an equivalent condition to the group configuration is expressed by a $4$-ary relation with additional properties. This uses ideas from \cite{CPS}.
\item As discussed in Section \ref{type definable}, one could replace the type-definable homogeneous space $(G,\Sigma)$ by {\em a definable} local group $G_0\supseteq G$, together with a definable family of injections on definable sets containing $x_2$.
\item Also in Section \ref{type definable} a related question is considered: Can every definably connected type definable group $G$ be definably embedded into a definable group, possibly of the same dimension? This is true in the stable setting, but the stable proof fails in the o-minimal setting. A  positive answer is given when $G$ is either $1$-dimensional or has a simple Lie algebra.
\end{enumerate}
\end{rem}

\vspace{.2cm}
\noindent {\bf References and acknowledgements.}
\vspace{.1cm}
While this o-minimal point of view  has not been published before, the main ideas are all taken from the existing accounts of various forms of the group configuration. Here are some of the sources that I have used.

The Group Configuration Theorem in the stable setting and its model theoretic proof are due to Hrushovski, \cite{Hr} and \cite{Hr4}.  However, already in  \cite[Proposition 5.1]{Z2} Zil'ber saw, using very different arguments, that a certain $6$-points configuration in the strongly minimal countably categorical setting (see \cite[Lemma 3.3]{Z1}) leads to the presence of a definable group.

Accounts of Hrushovski's theorem in the superstable setting appear in  Bouscaren's article \cite{BousGpCf}, and in the stable setting in Pillay's book \cite[Theorem 5.4.5]{Pi}. For the abelian version of the configuration I also used \cite{Bays-Hils-Moosa}. The current description was influenced by discussions with Chernikov and Starchenko during our work on \cite{CPS}, and also with Eleftheriou and Hasson during work on \cite{EHP} and elsewhere. The idea of using infinitesimal types in the o-minimal setting goes back to the early 1990's and was inspired by preliminary versions of Zil'ber's book, \cite{Zbook}. Some ideas here have already appeared in \cite{group-interval}.

\section{O-minimal preliminaries }
For basic reference on o-minimality, see \cite{vdDries} and also \cite{mynotes}.

We work in an o-minimal structure $\CM=\la M,<,\ldots\ra $ which is assumed to be $|\kappa|^+$ saturated, where $\kappa=|L|$ is the cardinality of the language.  Since $M$ is linearly ordered it can be  equipped with the order topology, and $M^n$ with the product topology. All topological references below are to this topology.

The letters $A,B,C,D$ are reserved for subsets of $M$ of cardinality at most $\kappa$.

\subsection{Geometric dimension of subsets of $M$}

 As  already noted, the existence of a definable linear ordering implies that $acl(A)=dcl(A)$ for every $A\sub M$ and  $dcl(-)$ gives rise to dimension of (small) subsets of $M$ over a (small) set of parameters. Namely,
$$\dim(A/B)=\{|A_0|:A_0\sub A \mbox{ is maximal } dcl\mbox{-independent over $B$}\}.$$
(a set $A\sub M$ is {\em $dcl$-independent over $B$} if for every $a\in A$, $a\notin dcl((A\setminus \{a\})\cup B))$.

\vspace{.3cm}
\noindent{\bf The dimension formula}
\vspace{.2cm}

For tuples $a,b$ and $A\sub M$,
$$\dim(a,b/A)=\dim(a/Ab)+\dim(b/A).$$

Two tuples  $a$ and $b$ are {\em inter-definable over $A$} if $dcl(aA)=dcl(bA)$. By the dimension formula, it follows that $\dim(a/A)=\dim(b/A)$.  The set $A$ is said to be {\em independent from $C$ over $B$} if $\dim(A/B\cup C)=\dim(A/B)$. By the Exchange Principle it follows that $\dim (C/A\cup B)=\dim(C/B)$. From now one we write $BC$ instead of $B\cup C$.

\vspace{.3cm}

\subsection{Definable sets}

If $X\sub M^n$ is a definable set over $B$ then
the dimension of  $X$ can be defined by
$$\dim (X)=\max\{\dim(a_1,\ldots,a_n/B):(a_1,\ldots,a_n)\in X\}.$$

An equivalent definition for $\dim(X)$ is the maximal $k\leq n$ such that some projection of $X$ onto $k$ of its coordinates contains an open subset of $M^k$. If $a\in X$ and $\dim(a/B)=\dim(X)$ then say that \emph{$a$ is generic in $X$ over $B$}. The following fact is very useful here:
\begin{fact} \label{neighborhood} Assume that  $a\in M^n$ and  $A\sub B\sub M$. For every open neighborhood $U$ of $a$, there exists $C\supseteq A$, independent from $aB$ over $A$,   and a $C$-definable open $W\sub U$ containing $a$. In particular, $\dim(a/A)=\dim(a/C)$ and $\dim(aB/C)=\dim(aB/A)$.

\end{fact}
\proof Consider a rectangular open box $W\sub U$ containing $a$ of the form $W=\Pi_{i=1}^n (a_i,b_i)$, with $\dim(a_1,\ldots,a_n,b_1,\ldots,b_n/aB)=2n$ (it is not hard to see, using saturation, that indeed such a box exists).
It follows that $a_1,\ldots,a_n,b_1,\ldots,b_n$ is independent from $aB$ over $A$. Let $C=Aa_1\cdots,a_nb_1\cdots b_n$.\qed

The following is \cite[Corollary 4.1.9]{vdDries}.
\begin{fact}\label{interior} Let $X\sub Y\sub M^n$ be definable sets and endow $Y$ with the subspace topology. Then $\dim(X\setminus Int_Y(X))<\dim Y$ (here $Int_Y(-)$ is the interior with respect to the topological space $Y$). In particular, if $a$ is generic in $Y$ over $A$ (assuming that $Y$ is $A$-definable) then for every $A$-definable $X\sub Y$ which contains $a$, we have $a\in Int_Y(X)$.

 \end{fact}

\subsection{Infinitesimal neighborhoods and infinitesimal loci}

The structure  $\CM$ is still assumed to be a $\kappa^+$-saturated structure.

\begin{defn} For $a\in M^n$,
{\em the $\CM$-infinitesimal neighborhood of $a$} is the partial type over $M$ consisting of all $M$-definable open subset of $M^n$ which contain $a$. For $\mathcal N$ an $|M|^+$-saturated elementary extension of $\CM$, it is often
convenient to identify this type with its realization in $\CN$.  Suppressing the dependence on $\CM$ the type is denoted by $\mu_a$ and its realization by $\mu_a(\CN)$. Since open boxes form a basis for the product topology, it follows that  if $a=(a_1,\ldots,a_n)$ then $\mu_a=\mu_{a_1}\times \cdots\times \mu_{a_n}$.

If $X\sub M^n$ is an $\CM$-definable set and $a\in X$ then let $\mu_a(X)$ denote the set $\mu_a(\CN)\cap X(\CN)$ (namely, the partial type $\mu_a\cup\{\varphi(\bar x)\}$, where $\varphi$ defines $X$).
\end{defn}

\begin{fact} \label{mu1} Assume that  $Q$ is $A$-definable and  $a$ is generic in $Q$ over $A$. Then for every $A$-definable set $R$, if $a\in R$ then $\mu_a(Q)\sub \mu_a(R)$. In particular, if $\dim Q=\dim R$ then $\mu_a(Q)=\mu_a(R)$.
\end{fact}
\proof First consider $Q\cap R$, which is also $A$-definable. By Fact \ref{interior}, $a\in Int_Q(Q\cap R)$ so there exists an $A$-definable open $U\ni a$ such that $U\cap Q\sub R\cap Q$. It follows that $U\cap R\cap Q=U\cap Q$ and therefore $\mu_a(R\cap Q)=\mu_a(Q)$.  Hence, $\mu_a(Q)\sub \mu_a(R)$. If $\dim Q=\dim R$ then $a$ is also generic in $R$ over $A$, hence  $\mu_a(R)\sub\mu_a(Q)$.\qed

\begin{defn} Given $a\in M^n$, $A\sub M$, and $Q\sub M^n$ an $A$-definable set such that $a$ is generic in $Q$ over $A$, one calls $\mu_a(Q)$ {\em the infinitesimal locus of $a$ over $A$, with respect to $\CM$}, and denote it by $\mu^{\CM}(a/A)$. If the ambient structure $\CM$ is clear then the notational reference to $\CM$ is omitted.
\end{defn}

\noindent {\bf Remarks} \begin{enumerate}

\item  By  Fact \ref{mu1}, the notion of $\Loc{a}{A}$ does not depend on the choice of $Q$.
\item  If $A\sub B \sub M$ and $B$ is independent from $a$ over $A$ then $\Loc{a}{A}=\Loc{a}{B}.$
\item If $a$ and $b$ are independent over $A$ then $\Loc{a,b}{A}=\Loc{a}{A}\times \Loc{b}{A}.$
\end{enumerate}

Here are some examples of infinitesimal loci.
\begin{example} Let $\CM=\la R,<,+,\cdot\ra$ be a $\aleph_1$-saturated real closed field, $R\supseteq \mathbb R$.
\begin{enumerate}

\item For $a\in R\setminus \mathbb R$,  $\Loc{a}{\mathbb R}=\Loc{a}{\emptyset}=\mu_a$ is the $\CM$-infinitesimal neighborhood of $a$,  namely the intersection of all open intervals $(a-\epsilon,a+\epsilon)$, $\epsilon\in R$.
\item If $a\in \mathbb R$ then $\Loc{a}{\mathbb R}$ is the partial type $\mu_a\cup\{x=a\}$ which is equivalent to the isolated type $x=a$.
\item If $a=(b,b^2)\in R^2$ with $b\notin \mathbb R$ then $\Loc{a}{\mathbb R}$ is $\mu_a\cup\{y=x^2\}$, the intersection
of $\mu_a$ with the parabola $y=x^2$. .
\item If $a=(a_1,a_2)\in R^2$ with $\dim(a/\mathbb R)=2$ then $\Loc{a}{A}=\mu_a$ the infinitesimal neighborhood of $a$ in $R^2$.
\end{enumerate}
\end{example}

 By Fact \ref{mu1},

\begin{fact}\label{loc-type} For every $a\in M^n$,
 $$\Loc{a}{A}\vdash tp(a/A).$$
 \end{fact}

 By the saturation of $\CM$, the types $\Loc{a}{A}$ and $tp(a/A)$ are logically equivalent
 if and only if $a \in dcl(A)$.

A  function $f:\Loc{a}{A}\to \Loc{b}{A}$ is called  {\em $A$-definable} if there exists an $A$-definable set $X$ such that the graph of $f$ equals $X\cap \Loc{a}{A}\times \Loc{b}{A}$.
Using compactness, it is not hard to see that in this case $X$ itself can be chosen to be the graph of a function.
\begin{lemma} \label{important} Let $a\in M^n$, $b\in M^k$.
\begin{enumerate}
\item The set $\Loc{ab}{A}$ is a subset of $\Loc{a}{A}\times \Loc{b}{A}$ and projects onto both.
\item  If $b\in dcl(aA)$  then the set $\Loc{ab}{A}$ is the graph of an $A$-definable  surjective function from $\Loc{a}{A}$ onto $\Loc{b}{A}$.

\item If $a$ and $b$ are inter-definable over $A$ then
$\Loc{ab}{A}$ is the graph of an $A$-definable bijection between $\Loc{a}{A}$ and $\Loc{b}{A}$.

\item Assume that  $a,b$ are independent over $A$ and $b,c$ are inter-definable over $aA$.
Then $\Loc{abc}{A}$ is the graph of an $A$-definable $F:\Loc{a}{A}\times \Loc{b}{A}\longrightarrow \Loc{c}{A}$ and  for every $a'\in \Loc{a}{A}$ , the map $F(a',-)$ is a bijection of $\Loc{b}{A}$ and $\Loc{c}{A}$.
\end{enumerate}
\end{lemma}

\proof (1) Assume that $X\ni (a,b)$ is an $A$-definable set such that $(a,b)$ is generic in $X$ over $A$.
We want to show that for every $b'\in \Loc{b}{A}$ there exists $a'\in \Loc{a}{A}$ with $(a',b')\in X$.

Let $U\ni a$ be an $M$-definable open $U\ni a$. By Fact \ref{neighborhood}, there exists an open $V\sub U$ containing $a$, defined over a parameter set $C$ independent
from $a,b$ over $A$. Hence, $\Loc{ab}{C}=\Loc{ab}{A}$

We have $(a,b)\in X$ and $a\in V$ and by Fact \ref{loc-type}, $\Loc{b}{C}\vdash tp(b/C)$. Hence, for every
$b'\in \Loc{b}{C}$ we have $\exists x\in V \,(x,b')\in X$. This is true for every $M$-definable open $U\ni a$
hence by compactness for every $b'\in \Loc{b}{C}$ there exists $a'\in \Loc{a}{A}$ with $(a',b')\in X$. This shows that the set $X\cap \Loc{a}{A}\times \Loc{b}{A}$,
which by definition equals to $\Loc{ab}{A}$, projects onto $\Loc{b}{A}$. By symmetry  it also projects surjectively onto $\Loc{a}{A}$.

(2) We repeat the same argument but since $b\in dcl(aA)$ one may also choose $X$ to be the graph of a function.

To obtain (3), we apply (2) to both coordinates.

(4) By (2), the set $\Loc{abc}{A}$ is the graph of a surjective function from
 $\Loc{ab}{A}=\Loc{a}{A}\times\Loc{b}{A}$ (by the independence of $a$ and $b$)
onto $\Loc{c}{A}$. For the same reason this set is  also the graph of a surjective function from $\Loc{a}{A}\times\Loc{c}{A}$ onto $\Loc{b}{A}$. The conclusion now follows.\qed

\begin{rem}  Some of the statements in Section 3.6.4 of \cite{Zbook}  resemble the above Lemma \ref{important}, despite the fact that both settings are very different.
\end{rem}

\subsection{Germs of definable functions}
Going back to Hrushovski, all proofs of the group configuration theorem in the stable case make use   of ``germs'' of  definable functions. The term itself originates in the theory of continuous or analytic functions, where local equality refers  to some ambient topology. O-minimality allows us to return to this original meaning.

\begin{defn} Given a definable set $X\sub M^n$ and definable maps $f,g:X\to M^k$, say that \emph{$f$ and $g$ have the same germ at $x_0\in X$}, $f\sim_{x_0} g$, if there exists an open $U\ni x_0$ such that for every $x\in U\cap X$, $f(x)=g(x)$.
\end{defn}

\begin{lemma} \label{germs}Let $x_0$ be a generic point in $X\sub M^n$ over $A$ and $\{f_t:t\in T\}$  an $A$-definable family of maps from $X$ into $M^k$. If $t_0\in T$ is generic and independent from $x_0$ over $A$ then there exist  definable open sets $W_0\ni t_0$, $U_0\ni x_0$, such that for every $t_1,t_2\in W_0\cap T$, and $x\in U_0$, if $f_{t_1}\sim_{x_0}f_{t_2}$ then $f_{t_1}| U_0\cap X=f_{t_2}| U_0\cap X$. Moreover, the sets $W_0$ and $U_0$  can be chosen to be defined over  parameters  independent from $x_0,t_0$ over $A$,

In particular, for every $t_1,t_2\in \Loc{t_0}{A}$, $$f_{t_1}\sim_{x_0}f_{t_2} \Leftrightarrow f_{t_1}| \Loc{x_0}{A}=f_{t_2}| \Loc{x_0}{A}\Leftrightarrow f_{t_1}|U_0=f_{t_2}|U_0.$$
\end{lemma}
\proof Without loss of generality, $A=\emptyset$.  To simplify notation write $t_1\sim_{x_0} t_2$ instead of $f_{t_1}\sim_{x_0}f_{t_2}$, and let $[t]_{x_0}$ denote the equivalence class.

 Assume that $\dim[t_0]_{x_0}=r$ and choose $t_1$ generic in $[t_0]_{x_0}$ over $t_0, x_0$. Let us see that $t_0$ is generic in $[t_0]_{x_0}=[t_1]_{x_0}$ over $t_1x_0$.

 Indeed,
$$\dim(t_0/t_1x_0)+\dim (t_1/x_0)=\dim(t_0t_1/x_0)=\dim(t_1/t_0x_0)+\dim(t_0/x_0)=r+\dim(T).$$
Since $\dim(t_1/x_0)\leq \dim T$, then $\dim(t_0/t_1x_0)\geq r$, but since $t_0\in [t_1]_{x_0}$
then  $\dim(t_0/t_1x_0)=r$, as claimed.

Let $U\ni x_0$ be an open set defined over additional parameters $B$ which are independent from $t_0,t_1,x_0$, such that $f_{t_0}|U=f_{t_1}|U$. Since $t_0$ is generic in $[t_1]_{x_0}$ over $Bx_0t_1$ there exists a $B$-definable open $W\ni t_0$ such that for every $t'\in W\cap T$, if $t'\sim_{x_0} t_1$ then $f_{t'}|U=f_{t_1}|U$.  But then, also
\begin{equation}\label{formula} \forall t'\in W\cap T\, (t'\sim_{x_0} t_0\longrightarrow f_{t'}|U=f_{t_0}|U).\end{equation}

Let $\phi(x_0,t_0)$ be a formula over $B$ expressing (\ref{formula}). Since $x_0$ and $t_0$ are independent over $B$, there are neighbrohoods $W_0\ni t_0$ and $U_0\ni x_0$ such that every $t_1\in W_0\cap T$ and $x_1\in U_0\cap X$ satisfy $\phi(x_1,t_1)$. Finally, one may replace $W_0$ and $U_0$ by $W_0\cap W$ and $U_0\cap U$, respectively, to obtain the required result.\qed

 \section{The configuration theorem}
 Assume that $\CM$ is an $|\CL|^+$-saturated structure which eliminates imaginaries.
\begin{defn}\label{defn1} {\em An $(m,k)$- homogenous space configuration in $\CM$ over $A\sub M$}, or just {\em an $(m,k)$-configuration over $A$},  is $6$-tuple of  tuples $(a_1,a_2,a_3,x_1,x_2,x_3)$ from $\CM$, such that (we allow redundancies in the clauses below)
\begin{enumerate}[(i)]
\item For $i=1,2,3$, $\dim(a_i/A)=m$ and  $\dim(x_i/A)=k$.
\item $\dim(x_1,x_2,x_3/A)=3k$ and for $i\neq j$, $\dim(a_i,a_j/A)=\dim(a_1,a_2,a_3/A)\!=2m$.
\item For distinct $i,j,k$, $\dim(x_i,x_j/a_kA)=k$.
\item For $i\neq j$, and any $k$,  $\dim(a_i,a_j,x_k/A)= 2m+k$ and
$\dim(a_i, x_i,x_j/A)=m+2k$.
\end{enumerate}

The configuration is called {\em minimal} if for every  $B\supseteq A$ independent from $V$  over $A$, and for every  $b_1,b_2,b_3$, with $b_i\in dcl(Ba_i)$, if $(b_1,b_2,b_3,x_1,_2,x_3)$ is an $(m',k)$-configuration over $B$ (so in particular, $m'\leq m$) then $m'=m$.
\end{defn}

 It is common to read-off the above information from the diagram below as follows: In addition to (i), every triples which is not co-linear is independent over $A$ and also, whenever $x_i,x_j, a_k$, $i\neq j$,  is  a co-linear triple then $x_i\in dcl(x_j,a_k,A)$, and  each element of the co-linear triple $a_1,a_2,a_3$ is in the definable closure over $A$ of the other two.
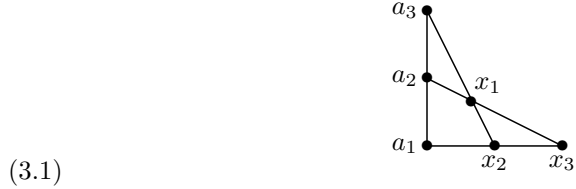
\begin{equation}
  \begin{tikzpicture}[scale=0.9]
\draw[black,line  width=0.2mm] (0,0) -- (0,2);
\draw[black,line  width=0.2mm] (0,0) -- (2,0);
\draw[black,line  width=0.2mm] (1,0) -- (0,2);
\draw[black,line  width=0.2mm] (0,1) -- (2,0);
\node at (0,0) {$\bullet$};
\node[left,scale=1.0] at (0,0) {$a_1$};
\node at (1,0) {$\bullet$};
\node[below,scale=1.0] at (1,0) {$x_2$};
\node at (2,0) {$\bullet$};
\node[below,scale=1.0] at (2,0) {$x_3$};
\node at (0,1) {$\bullet$};
\node[left,scale=1.0] at (0,1) {$a_2$};
\node at (0,2) {$\bullet$};
\node[left,scale=1.0] at (0,2) {$a_3$};
\node at (0.65,0.65) {$\bullet$};
\node at (0.9,0.9) {$x_1$};
\end{tikzpicture}
\end{equation}

\begin{rem} It follows from the dimension formula that whenever $V$ is an $(m,k)$-configuration over $A$ then $m\geq k$.
\end{rem}

The canonical example of an $(m,k)$-configuration is given by a group action as follows:
\begin{example}\label{gpexample}
Let $G$ be an type-definable group acting definably and transitively on a definable $k$-dimensional set $X$.
Let $g_1,g_2$ be independent and generic elements in $G$ and $x\in X$ generic and independent over $a_1,a_2$.
It is not hard to verify that $(g_1,g_2,g_2g_1,x,g_1\cdot x, g_2g_1\cdot x)$ is an $(m,k)$-configuration. Namely,
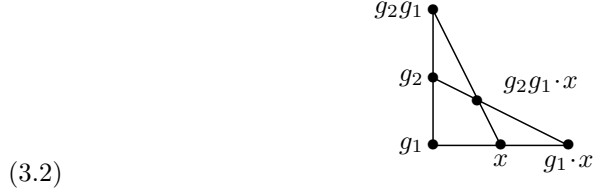
\begin{equation}
  \label{eq:2}
  \begin{tikzpicture}[scale=0.9]
\draw[black,line  width=0.2mm] (0,0) -- (0,2);
\draw[black,line  width=0.2mm] (0,0) -- (2,0);
\draw[black,line  width=0.2mm] (1,0) -- (0,2);
\draw[black,line  width=0.2mm] (0,1) -- (2,0);
\node at (0,0) {$\bullet$};
\node[left,scale=1.0] at (0,0) {$g_1$};
\node at (1,0) {$\bullet$};
\node[below,scale=1.0] at (1,0) {$x$};
\node at (2,0) {$\bullet$};
\node[below,scale=1.0] at (2,0) {$g_1\!\cdot\! x$};
\node at (0,1) {$\bullet$};
\node[left,scale=1.0] at (0,1) {$g_2$};
\node at (0,2) {$\bullet$};
\node[left,scale=1.0] at (0,2) {$g_2g_1$};
\node at (0.65,0.65) {$\bullet$};
\node[right,scale=1] at (0.9,0.9) {$g_2g_1\!\cdot\! x$};
\end{tikzpicture}
\end{equation}
If the action of $G$  is faithful then the configuration is minimal.

Note that the parameter set over which the configuration is defined depends only on the
parameters defining the group operation of $G$ and not on those defining the universe of $G$. This will be used below.
\end{example}

 Here is  the o-minimal group-configuration theorem.
  \begin{theorem} \label{main-thm}Assume that $V=(a_1,a_2,a_3,x_1,x_2,x_3)$ is an $(m,k)$-configuration over $A$. Then there exists  an $m'$-dimensional type-definable group $(G,\star)$, with $k\leq m'\leq m$,
acting definably and transitively on some $k$-dimensional partial type $\Sigma$ over $M$.

Moreover, there is a parameter set $B\supset A$ independent from $V$ over $A$, and there are elements $g_1,g_2,g_3\in G$, and $y_1,y_2,y_3 \models \Sigma$ such that for each $i=\{1,2,3\}$, we have $g_i\in dcl(a_iB)$,  and $dcl(x_iB)=dcl(y_iB)$, and in addition $(g_1,g_2,g_3,y_1,y_2,y_3)$ is an $(m',k)$-configuration over $B$ inside $G$ as in Example \ref{gpexample}.

If, in addition, the configuration is minimal then $m'=m$ and each $g_i$ is inter-definable with $a_i$ over $B$.
\end{theorem}

\begin{rem} Although the universe of the type-definable group $G$ is definable over $M$ (in fact, it equals $\Loc{a_1}{A}$), the group operation will be defined over $B$ and therefore the configuration
coming from Example \ref{gpexample} is a configuration over $B$.
\end{rem}
\noindent {\em Proof of the theorem.}  Note that by assumptions, for the co-linear triple $\{a_1,a_2,a_3\}$, each two of the $a_i$'s are inter-definable over the third one, and for the other co-linear triples, the $x_i$'s are inter-definable over the $a_j$.
Here is the main idea of the proof:
By applying Lemma \ref{important} (3),  one obtains  $A$-definable functions $F,L,H,K$ satisfying:
$$F(a_2,a_1) =a_3\,\,\, ,\,\,\, L(a_1,x_2)=x_3\mbox{ and }  H(a_3,x_2)=x_1=K(a_2,x_3).$$

Putting together the above functional equation, one  obtains the following associativity-like law.

\begin{equation}\label{assoc}
\begin{array}{c}
  H(F(a_2,a_1),x_2)=K(a_2,L(a_1,x_2)), \mbox{ with }\\ \\ \dim(a_i/A)=m, \mbox{ for } i=\!1,\!2, \,\, \dim(x_2/A)=k
  \mbox{ and } \dim(a_1,a_2,x_2/A)=2m+k\\
\end{array}
 \end{equation}






The existence of definable functions $H,F.K.L$ and points $a_1,a_2,x_2$ satisfying (\ref{assoc}), could be the appropriate version of the group configuration in the o-minimal setting, as it reflects  the functional nature  of o-minimal objects.

Let us now go into the details.

$(i) $ \hspace{.2cm}
First apply Fact \ref{important} (4)  to $a_1,a_2,a_3$ and obtain an  $A$-definable function $F(-,-)$, satisfying
$$F:\Loc{a_2}{A}\times \Loc{a_1}{A} \rightarrow \Loc{a_3}{A} \mbox { and } F(a_2,a_1)=a_3.$$ Moreover, for every $a_1'\in \Loc{a_1}{A}$, $F(-,a_1')$ is a bijection. The analogous fact holds for $F(a_2',-) $.

$(ii)$  \hspace{.2cm} Next, apply the same result to the triple $a_1, x_2,x_3$ and obtain an $A$-definable $L(-,-)$ satisfying
 $$L:\Loc{a_1}{A}\times \Loc{x_2}{A}\rightarrow \Loc{x_3}{A} \mbox { and } L(a_1,x_2)=x_3,$$   such that for every $a_1'\in \Loc{a_1}{A}$, the function $L(a_1',-)$ is a bijection.

$(iii)$ \hspace{.2cm} We also have $$H:\Loc{a_3}{A}\times \Loc{x_2}{A} \rightarrow \Loc{x_1}{A} \mbox{ and } H(a_3,x_2)=x_1,$$  such that for every $a_3'\in \Loc{a_3}{A}$ the function $H(a_3',-) $ is a bijection.

$(iv)$ \hspace{.2cm} Finally, we have $$K:\Loc{a_2}{A}\times \Loc{x_3}{A}\rightarrow \Loc{x_1}{A} \mbox{ and } K(a_2,x_3)=x_1,$$ such that for every $a_2'\in \Loc{a_2}{A}$, $K(a_2',-)$ is a bijection.
We conclude
\begin{equation} \label{functional eq} H(F(a_2,a_1), x_2)=K(a_2,L(a_1,x_2)) \end{equation}

\subsection{Composition of families of functions}
\newcommand{\CH}{\mathcal H}
\newcommand{\CG}{\mathcal G}
\newcommand{\CK}{\mathcal K}

For $a_3'\in \Loc{a_3}{A}$,  let $$h_{a_3'}(-)=H(a_3',-):\Loc{x_2}{A}\to \Loc{x_1}{A}.$$ Similarly, write $$k_{a_2'}:=K(a_2',-):\Loc{x_3}{A}\to \Loc{x_1}{A} \mbox{ and }\ell_{a_1'}:=L(a_1',-):\Loc{x_2}{A}\to \Loc{x_3}{A}$$ and $$ f_{a_2'}:=F(a_2,-):\Loc{a_1 }{A}\rightarrow \Loc{a_3}{A}.$$

Consider now the three families of functions:
 $$\CH=\{h_{a_3'}:a_3'\in \Loc{a_3}{A}\}\, ;\, \mathcal K=\{k_{a_2'}:a_2'\in \Loc{a_2}{A}\}\, ; \, \mathcal L=\{\ell_{a_1'}:a_1'\in \Loc{a_1}{A}\} .$$

Using (\ref{functional eq}),
\begin{claim}\label{comp1.1}
 (i) For every $k\in \CK$ and $\ell\in \CL$, there exists $h\in \CH$ such that $k\circ \ell=h$.
\\

(ii) For every $h=h_{a_3'}\in \CH$, and $k=k_{a_2'}\in \CK$ there exists $\ell=\ell_{f_{a_2}^{-1}(a_3')}\in \CL$
such that $k\circ \ell=h$.
\\

(iii)  For every $\ell\in \CL$ and $h\in \CH$ there exists $k\in \CK$ such that $k\circ \ell=h$.
\end{claim}
\proof Let us see for example (i). If $k=k_{a_2'}$ and $\ell=\ell_{a_1'}$ then $h=h_{F(a_2',a_1')}$.\qed

Note that the above claim implies that the families $\CK$, $\CL$ and $\CH$ all have the same dimension (at most $k$), although $\CH$ is obtained as the composition of functions from $\CK$ and $\CL$.  This is another indication for the presence of a group.

From here on composition $f\circ g$  is written as $fg$.
We are now ready to define the intended group:
$$G=\{\ell_1^{-1}\ell_2:\ell_i\in \mathcal L\},$$ a family of permutations of $\Loc{x_2}{A}$.

\begin{prop}\label{group} The set $G$ is a group with respect to composition.\end{prop}

\proof The main tool towards the proof is:
\begin{lemma} \label{presentation} For every $\ell_1,\ell_2,\ell_3\in \mathcal L$ there exists $\ell_4\in \CL$ such that $\ell_1^{-1}\ell_2=\ell_3^{-1} \ell_4$.\end{lemma}
\proof By Claim \ref{comp1.1} (i), we can write $\ell_3=k_3^{-1}h_3$ for some $k_3\in \CK$, $h_3\in \CH$. Then write $\ell_1=k_1^{-1}h_3$ for some $k_1\in \CK$, and finally $\ell_2=k_2^{-1}h_2$ for some $h_2\in \CH$. Thus
$\ell_3 \ell_1^{-1}\ell_2=k_3^{-1} h_2 $ which by Claim \ref{comp1.1} (ii), equals some $\ell_4\in \mathcal L$, thus proving the lemma.\qed

To see that $G$ is indeed a group one needs to see that for every $\ell_1,\ell_2,\ell_3,\ell_4$, the function $\ell=\ell_1^{-1}\ell_2\ell_3^{-1}\ell_4$ is in $G$.
By the above lemma, $\ell_3^{-1} \ell_4$ can be written as $\ell_2^{-1}\ell_5$ for some $\ell_5\in G$ and then $\ell=\ell_1^{-1}\ell_5$ is in $G$, thus proving the proposition.\qed

 By its definition the group $G$ acts on $\Loc{x_2}{A}$, as $(\ell_1^{-1}\ell_2)\cdot x=\ell_1^{-1}(\ell_2(x))$.
\begin{lemma} The action of $G$ on $\Loc{x_2}{A}$ is transitive, so in particular $\dim G\geq k$.\end{lemma}
\proof
 The first  claim is that for every $(x_2',x_3')\in \Loc{x_2}{A}\times \Loc{x_3}{A}$, there exists $a_1'\in\Loc{a_1}{A}$ such that $\ell_{a_1'}(x_2')=(x_3').$ Indeed, by Lemma \ref{important}, $\Loc{a_1x_2x_3}{A}$, which is exactly $Graph(L)\cap \Loc{a_1}{A}\times \Loc{x_2}{A}\times \Loc{x_3}{A}$,  projects onto $\Loc{x_2}{A}\times \Loc{x_3}{A}$. The claim follows.

Given $(x_2',x_2'')\in \Loc{x_2}{A}\times \Loc{x_2}{A}$, we can therefore find $\ell_1,\ell_2 \in \CL$ such that $\ell_1(x_2'')=x_3$ and $\ell_2(x_2')=x_3$. It follows that $\ell_1^{-1}\ell_2(x_2')=x_2'',$ hence $G$ acts transitively.\qed

The next goal is to identify $G$ as a type-definable group. First fix $a\in \Loc{a_1}{A}$ independent form all parameters
mentioned thus far (including $V$ itself). By Lemma \ref{presentation}, $G=\{\ell_{a}^{-1}\ell_{a_1'}:a_1'\in \Loc{a_1}{A}\}$.

  By Lemma \ref{germs}, there exist a definable   open $U_0\ni x_2$ and a definable open $W\ni a_1$,  possibly defined over additional independent parameters, such that for every $a'_1,a_1''\in W$, $\ell_{a_1'}=\ell_{a_1''}$ (as functions on $\Loc{x_2}{A}$) iff for every $x\in U_0$, $L(a_1',x)=L(a_1'',x)$. We may also assume that $\ell_a^{-1}|U_0$ is injective.
  
It follows that there exists a definable equivalence relation $E$, defined over independent parameters $C$, such that for every $a_1',a_1''\in W$,
$\ell_a^{-1}\ell_{a_1'}=\ell_a^{-1}\ell_{a_1''}$ if and only if $a_1'E a_2'$. Denote by $[a_1']$ the equivalence class of $a_1'$ with respect to $E$.
Because $\CM$ eliminates imaginaries there is a definable bijection $\alpha$ between the quotient space $W/E$ and a definable set $Z$.  We may choose $W$ small enough so that the dimension of $W/E$ is exactly $m'=\dim G=\dim Z$.

 The element $\alpha(a_1)$ is generic in $Z$ over $C$ and  as in the proof of Lemma \ref{important}(1), $\alpha(\Loc{a_1}{A})=\alpha(\Loc{a_1}{C})=\Loc{\alpha(a_1)}{C}$, so the map $\alpha$ induces a bijection between elements of $G$ which are functions $\ell_a^{-1}\ell_{a_1'}$ and  $\alpha(a_1')\in \Loc{\alpha(a_1)}{C}.$ This shows that $G$ can indeed be realized as a type-definable set $\Loc{\alpha(a_1)}{C}$.

Finally,  one needs to find a corresponding configuration in $G$. For that, it is still convenient to view $G$ as the set $\ell_a^{-1}\ell_{a_1'}$, for $a_1'\in \Loc{a_1}{A}$. Fix $b\in \Loc{a_3}{A}$ generic  over $A$ and all parameters thus far (including the configuration $V$ and $a$).


Define the elements of the $G$-configuration as follows:
$$g_1:=\ell_a^{-1}\ell_{a_1}\, \, g_2:=h_b^{-1}k_{a_2}\ell_a\, ,\,  g_3:=g_2g_1=h_b^{-1}k_{a_2}\ell_{a_1}.$$
We have $g_1\in dcl(Aa a_1)$, $g_2\in dcl(Aba_2)$ and  by (\ref{functional eq}), $k_{a_2}\ell_{a_1}=h_{a_3}$, hence $g_3\in dcl(Aba_3)$.

Now let $$y_2:=x_2\, ;\,y_3:=g_1\cdot y_2=\ell_a^{-1}\ell_{a_1}(x_2)=\ell_a^{-1}(x_3).$$
and $$ y_1:=g_3\cdot y_2 h_b^{-1}h_{a_3}(y_2)=h_b^{-1}h_{a_3}(x_2)=h_b^{-1}(x_1).$$ Hence, $y_i\in dcl(Aabx_i)$,
and by its definition and the independence of $ab$ over all parameters,
$V_G=(g_1,g_2,g_3,y_1,y_2,y_3)$ is a group configuration as in Example \ref{gpexample}  over $Cb$.
This ends the proof of the main statement in Theorem \ref{main-thm}.

Assume now that the configuration $V$ is minimal. Then, the above configuration shows that necessarily $m'=m$. Indeed, if not then $V_G$ contradicts the minimality of $V$. Also, the equivalence class of $a_1$ (with respect to the germ relation $E$ above) is finite. It is not difficult to see that in fact, the equivalence class of each $a_1'\in \Loc{a_1}{A}$ contains a single element $a_1'$.
This ends the proof of Theorem \ref{main-thm}.\qed




\subsection{The abelian configuration}

\begin{defn} An {\em abelian $(m,k)$ configuration} is a configuration as in Definition \ref{defn1} with  an additional node $x_4$ such that $\dim(x_4/A)=k$,
and two more edges, connecting it to $a_2,x_2$ on one hand and to $a_1,x_3$ on the other. Namely,
$x_1, x_4$ are inter-definable over $a_1$ and $x_2,x_4$ are inter-definable over $a_2$. Equivalent way of describing it is by saying that $(a_1,a_2,a_3, x_2,x_4,x_1)$ is also a group configuration, with no further dependencies.
\end{defn}
\begin{equation}
  \begin{tikzpicture}[scale=0.9]
\draw[black,line  width=0.2mm] (0,0) -- (0,2);
\draw[black,line  width=0.2mm] (0,0) -- (2,0);
\draw[black,line  width=0.2mm] (1,0) -- (0,2);
\draw[black,line  width=0.2mm] (0,1) -- (2,0);
\draw[black,line  width=0.2mm] (0,1) -- (1,0);
\draw[black,line  width=0.2mm] (0,0) -- (.65,.65);
\node[below,scale=1.0] at (0.5,0.5) {$x_4$};
\node at (0.5,0.5){$\bullet$};
\node at (0,0) {$\bullet$};
\node[below,scale=1.0] at (0,0) {$a_1$};
\node at (1,0) {$\bullet$};
\node[below,scale=1.0] at (1,0) {$x_2$};
\node at (2,0) {$\bullet$};
\node[below,scale=1.0] at (2,0) {$x_3$};
\node at (0,1) {$\bullet$};
\node[left,scale=1.0] at (0,1) {$a_2$};
\node at (0,2) {$\bullet$};
\node[left,scale=1.0] at (0,2) {$a_3$};
\node at (0.65,0.65) {$\bullet$};
\node at (0.9,0.9) {$x_1$};
\end{tikzpicture}
\end{equation}

\begin{theorem} Given an abelian $(m,k)$-configuration, the group $G$ obtained
in Theorem \ref{main-thm} is abelian.
\end{theorem}
\proof In addition to the functional equations in (\ref{assoc}) there are now two additional definable functions, call them $R$ and $S$,  such that
$R(a_1,x_4)=x_1$ and $S(a_2,x_2)=x_4$, so $x_1=R(a_1,S(a_2,x_2))$. Together with the equations from Definition \ref{functional eq} we obtain:
\begin{equation}\label{abfunctional eq}  H(F(a_2,a_1), x_2)=R(a_1,S(a_2,x_2))=K(a_2,G(a_1,x_2)) .\end{equation} And translating it to compositional equations, we have
$$h_{F(a_2,a_1)}=r_{a_1}s_{a_2}=k_{a_2}\ell_{a_1}, \mbox{ as  functions from $\Loc{x_2}{A}$ to $\Loc{x_1}{A}$}, $$ and for every $a_1',a_2'$ in $\Loc{a_1}{A}$ and $\Loc{a_2}{A}$, respectively, we have
\begin{equation} \label{ab-eq}h_{F(a_2',a_1')}=k_{a_2'}\ell_{a_1'}=r_{a_1'}s_{a_2'}.\end{equation}

Let $\mathcal R=\{r_{a_1'}:a_1'\in \Loc{a_1}{A}\}$ and $\mathcal S=\{s_{a_2'}:a_2'\in \Loc{a_2}{A}\}$.
The main tool for the commutativity of $G$ is the following (taken from \cite{CPS}).

\begin{claim}\label{ab-claim} For any $h_1,h_2,h_3\in \CH$, $h_1h_2^{-1}h_3=h_3h_2^{-1}h_1$.\end{claim}
\proof By (\ref{ab-eq}), there are $a_1'\in \Loc{a_1}{A}$, $a_2'\in \Loc{a_2}{A}$ such that
$$h_3=r_{a_1'}s_{a_2'}=k_{a_2'}\ell_{a_1'}.$$

By Claim \ref{comp1.1}, and by (\ref{ab-eq}), there is $a_1''\in \Loc{a_1}{A}$, such that 
$$h_2=k_{a_2'}\ell_{a_1''}=r_{a_1''}s_{a_2'}.$$

Also, by Claim \ref{comp1.1}, and by (\ref{ab-eq}), there is $a_2''\in \Loc{a_2}{A}$, such that
$$h_1=k_{a_2''}\ell_{a_1''}=r_{a_1''}s_{a_2''}.$$

We now have $h_1h_2^{-1}h_3= k_{a_2''}\ell_{a_1'}$ and $h_3h_2^{-1}h_1= r_{a_1'}s_{a_2''}$, so by  (\ref{ab-eq}),
we have $h_1h_2^{-1}h_3=h_3h_2^{-1}h_1.$ \qed

We can now conclude that the group $G=\{\ell_1^{-1}\ell_2:\ell_i\in \CL\}$ is abelian. It is not hard to see, using (\ref{presentation}) that each element of $G$ can also be written as $h_1h_2^{-1}$, for $h_i\in \CH$, and as before for each $g\in G$ and $h\in \CH$ there exists $h_1\in \CH$ such that $g=h_1h^{-1}$.

Thus, in order to prove commutativity, one needs to show that for every $h_1,h_2,h_3\in \CH$,
$$h_1h_2^{-1}h_3h_2^{-1}=h_3h_2^{-1}h_1h_2^{-1}.$$ This is immediate from Claim \ref{ab-claim}.\qed

\section{$4$-ary relations and the group configuration}
\label{quadrangle}

While working on \cite{CPS}, a new way of viewing the data of the group configuration has emerged,
which is based on the relations on the $4$-tuple $(a_1,a_2,x_2,x_1)$
from Definition \ref{defn1}.

For simplicity, we assume here that our o-minimal structure $\CM$ eliminates imaginaries.

\begin{defn} For $k\geq 1$, a $4$-tuple $(a_1,a_2,a_3,a_4)$ is called {\em a $(k,k)$-quadrangle over $A$} if the following holds:
\begin{enumerate}
\item For each $i$, $\dim(a_i/A)=k$, each three of the $a_i$'s are independent over $A$ and $\dim(a_1a_2a_3a_4/A)=3k$.
\item Let $Q^*=\Loc{a_1a_2a_3a_4}{A}$. Then for every $(a_1',a_2',a_3',a_4'),(a_1'',a_2'',a_3'',a_4'')\in \Loc{a_1}{A}\times\Loc{a_2}{A}\times\Loc{a_3}{A}\times\Loc{a_4}{A}$, the following hold:

\vspace{.3cm}

\noindent (i) If $(a_1',a_2',a_3',a_4'),(a_1'',a_2'',a_3',a_4'), (a_1',a_2',a_3'',a_4'')\!\in \!Q^*$ then $ (a_1'',a_2'',a_3'',a_4'')\!\in\! Q^*$.

\noindent (ii) If $(a_1',a_2',a_3',a_4'),(a_1'',a_2',a_3',a_4''), (a_1',a_2'',a_3'',a_4')\in Q^*$ then $(a_1'',a_2'',a_3'',a_4'')\in Q^*$.

The quadrangle is called {\em abelian} if in addition:

\noindent (iii) If $(a_1',a_2',a_3',a_4'),(a_1',a_2'',a_3',a_4''), (a_1'',a_2',a_3'',a_4')\in Q^*$ then  $(a_1'',a_2'',a_3'',a_4'')\in Q^*$.

\end{enumerate}
\end{defn}

Notice that (1) implies that for every $(a_1',a_2')\in \Loc{a_1a_2}{A}$, $Q^*(a_1',a_2',-,-)$ is a bijection of $\Loc{a_3}{A}$ and $\Loc{a_4}{A}$ (see Lemma \ref{important}), call it $h^*_{a_1'a_2'}$ and (i) implies that if $h^*_{a_1'a_2'}$ and $h^*_{a_1''a_2''}$ agree at any point of $\Loc{a_3}{A}$ then they agree everywhere.

Similarly, for every  for every $(a_1',a_4')\in \Loc{a_1a_4}{A}$, $Q^*(a_1',-,-,a_4')$ is a bijection of $\Loc{a_2}{A}$ and $\Loc{a_3}{A}$, call it $\ell^*_{a_1'a_4'}$ and by (ii) it is determined by any point in its graph. Finally, there are functions $s_{a_1'a_3'}$ given by $Q^*(a_1',-,a_3',-)$, and if (iii) holds then they are determined by any points in their graph.

Assume now that $Q^*$ is a $(k,k)$-quadrangle over $A$. By compactness, there exists  $Q$  an $A$-definable set whose intersection with $\Loc{a_1}{A}\times\Loc{a_2}{A}\times\Loc{a_3}{A}\times\Loc{a_4}{A}$ equals $Q^*$, such that $Q$ is the graph of an injection whenever we fix two of its coordinates. Let $h_{a_1'a_2'}$ be the definable function $Q(a_1',a_2',-,-)$ and $\ell_{a_1'a_4'}$ be $Q(a_1',-,-,a_4')$.

Using Lemma \ref{germs}, there exists a open set $U_3\ni a_3$, definable over independent parameters, such that
$h_{a_1'a_2'}=h_{a_1''a_2''}$ iff the two functions agree on $U$. We let $[h_{a_1'a_2'}]$ denote the canonical parameter
of the definable function $h_{a_1'a_2'}|U_3$. Similarly, define the canonical parameters $[k_{a_1'a_4'}]$ and $[s_{a_1'a_3'}]$, using an an appropriately defined neighborhood of $a_2$. Assume that these two neighborhoods are defined over $C$.

Notice that if $(a_1,a_2,a_3,a_4)$ is a $(k,k)$-quadrangle then $[h_{a_1a_2}]\in dcl(a_1a_2C)\cap dcl(a_3a_4C)$
and $[k_{a_1a_4}]\in dcl (a_2a_3C)\cap dcl(a_1a_4C)$. If the quadrangle satisfies (iii) then also $[s_{a_1a_3}]\in dcl(a_1a_3C)\cap dcl(a_2a_4C)$.

 Without going into the details I state the following connection between $(k,k)$-configurations and $(k,k)$-quadrangles.
\begin{prop}
\begin{enumerate}
\item Assume that $(a_1,a_2,a_3,a_4)$ is a $(k,k)$-quadrangle over $A$. Then  the following is a $(k,k)$ configuration over $A$:
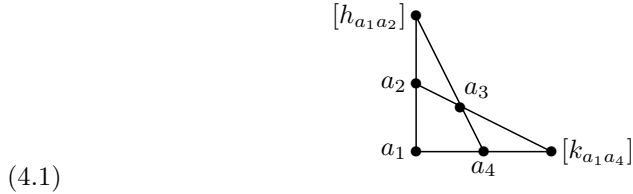
\begin{equation}
  \begin{tikzpicture}[scale=0.9]
\draw[black,line  width=0.2mm] (0,0) -- (0,2);
\draw[black,line  width=0.2mm] (0,0) -- (2,0);
\draw[black,line  width=0.2mm] (1,0) -- (0,2);
\draw[black,line  width=0.2mm] (0,1) -- (2,0);
\node at (0,0) {$\bullet$};
\node[left,scale=1.0] at (0,0) {$a_1$};
\node at (1,0) {$\bullet$};
\node[below,scale=1.0] at (1,0) {$a_4$};
\node at (2,0) {$\bullet$};
\node[right,scale=1.0] at (2,0) {$[k_{a_1a_4}]$};
\node at (0,1) {$\bullet$};
\node[left,scale=1.0] at (0,1) {$a_2$};
\node at (0,2) {$\bullet$};
\node[left,scale=1.0] at (0,2) {$[h_{a_1a_2}]$};
\node at (0.65,0.65) {$\bullet$};
\node at (0.9,0.9) {$a_3$};
\end{tikzpicture}
\end{equation}

If the quandrangle is abelian then by adding $[s_{a_1a_3}]$ one obtains an abelian configuration.

\item Assume that  $V=(a_1,a_2,a_3,x_1,x_2,x_3)$ is a $(k,k)$ configuration over $A$. Then
 $(a_1,a_2,x_1,x_2)$ is a $(k,k)$-quadrangle $A$. If $V$ can be extended to an abelian configuration then the quadrangle
is abelian.
\end{enumerate}
\end{prop}

\section{Type-definable groups}
\label{type definable}
\subsection{Local groups}

The notion of a ``a local group'' in the topological setting (not to be confused with ``a locally definable group'') replaces  the notion of ``a group-chunk'' from algebraic geometry. Without going into a formal definition, it is a topological space $X$ with a distinguished element $1\in X$, and a binary continuous operation $M$ defined in a neighborhood of $1$ such that $M(x,1)=M(1,x)=x$ and $M$ is associative when defined. In addition there is a unary continuous function in a neighborhood of $1$ which sends every element $g$ to an element $h$ such that $M(g,h)=M(h,g)=1$.

If $G$ is a type-definable group in an o-minimal structure then by the work of Marikova \cite{Marikova}, it admits a definable topology (namely a definable basis) making it into a topological group. Moreover, if $\CM$ expands a real closed field then by \cite[Theorem 3.7]{Marikova},  $G$ is definably isomorphic to a type-definable group in $M^n$ whose whose group topology is induced by $M^n$. Using logical compactness it is not hard to see that $G$ is contained in a definable local group, with respect to this topology.
Thus, in Theorem \ref{main-thm} one can replace the type-definable group $G$ with a definable local group, together with a definable family of injections which restrict to an action of $G$ on $\Loc{x_2}{A}$.

There are very few results, beyond Marikova's,  on general type-definable groups in o-minimal structures.   The work on Pillay's Conjecture required some analysis of certain type-definable subgroups of definable groups, and then in \cite{BaysPe} we examined the structure of $G^{00}$ for  a definably compact group $G$.

\subsection{Embedding type-definable groups into definable ones}
By  theorem of Hrushovski \cite[Theorem 5.18]{PoiGroups}, every type-definable group in a stable structure is contained
in a definable group. In addition, by Poizat, \cite[Theorem 5.17]{PoiGroups}, every type-definable subgroup of a stable group $H$ can be written as the intersection of definable subgroups of $H$. It follows from DCC
that every type-definable group in an $\omega$-stable strucutre is in fact definable.

What can be expected in the o-minimal setting? First note that the subgroup of infinitesimals $\mu_0=\{|x|<1/n:n\in \mathbb N\}$ of a real closed field $(R,+)$ is type-definable but cannot be written as  the intersection of definable subgroups. Hence, the analogous result to Poizat's fails. However, $\mu_0$ can still be definably embedded into the definably compact group $([0,1),+\, \mbox{mod}1)$. Thus, the following question arises. A type definable group $G$ is called {\em definably connected} if there is no definable set $U$ whose intersection with $G$ is a non-trivial clopen set in the sense of the group topology.

\vspace{.3cm}

\noindent  {\bf Question} {\em Can every definably connected type-definable group $G$ in an o-minimal structure  be definably embedded into a definable group $H$? Moreover, can one choose $H$ to be of the same dimension as $G$?}
\vspace{.3cm}

Below are some cases where a positive answer is given. Before stating the result recall that if $\CM$ expands a real closed field $R$ then every definable group $G$ has an associated Lie $R$-algebra $L(G)$ (see \cite{PPS1} for details). Since the definition of the Lie algebra uses only local information around the identity, the same construction works for type-definable groups, and thus if $G$ is type-definable then it has a corresponding Lie algebra $L(G)$ over $R$.
We say that $G$ is {\em Lie-simple} if $L(G)$ is a simple Lie algebra. When $G$ is definable then by \cite[Theorem 2.36]{PPS1}, $G$ is Lie-simple if and only if it is definably simple, namely has no definable normal subgroups.

\begin{prop} Let $G$ be a definably connected type-definable group in an o-minimal expansion of a real closed field, such that one of the two hold (i) $\dim G=1$ and $G$ is an ordered group,  or (ii) $G$ is Lie-simple.

 Then $G$ can be embedded into a definable group of the same dimension. Moreover, in case (ii) the definable group is semialgebraic.\end{prop}
\proof Assume $\dim G=1$. By assumption, $G$ is torsion-free. First note that $G$ is abelian. Indeed, for every $g\in G$, the centralizer $C_G(g)$ is an infinite, relatively definable subgroup of $G$. Because $\dim G=1$ it follows that $C_G(g)$ is open in $G$, and hence also closed, so by connectedness it must equal $G$. We write the group operation additively.

 As we noted above, we may assume that $G\sub M^n$ and its group topology is the $M^n$-topology. Since $G$ is definably connected it is given as the intersection of definably connected $1$-dimensional subsets $\bigcap I_i$ of $M^n$.  We may assume that $G$ is not definable itself and hence we may assume that the definable sets $I_i$ are (possibly after a definable bijection) all open intervals in $M$. Thus, we assume that the universe of $G$ is given as the intersection of open intervals $I_i$, and furthermore that the group operation makes each $I_i$ into a local group. Finally we may also assume that each $I_i$ is symmetric with respect to $+$ and we write it as $(-a_i,a_i)$.  By compactness, for each $I_i$, there exists $I_j\sub I_i$ such that $I_j+I_j=(-2a_j,2a_j)\sub I_i$. Fix such a pair $I_j\sub I_i$ and define on the interval $[-a_j,a_j)$ the following operation:
 
 \begin{equation}
 x\oplus y=\begin{cases}
 x+y & \mbox{ if } x+y\in [-a_j,a_j)\\
 x+y-2a_j & \mbox{ if } x+y\geq a_j\\
 x+y+2a_j & \mbox{ if } x+y<-a_j
 \end{cases}
 \end{equation}
 
 We obtain a definable (definably compact)  group $\hat G$ on the set $[-a_j,a_j)$, and by definition of $\oplus$, the group $G$ is a subgroup of  $\hat G$.
 This ends the proof of this case.



Assume now that $G$  is Lie-simple. Using  the adjoint embedding $Ad:G\to GL_m(R)$ (see\cite[Proof of Theorem 3.2]{PPS1})  one can view as  $G$ as a type-definable subgroup of $Aut(L(G))$ (see \cite[Claim 2.29]{PPS1}). By \cite[Claim 2.8]{PPS1}, $\dim(Aut(L(G))=\dim L(G)=\dim G$. The group $Aut(L(G))$ is clearly seialgebraic (one can explicitly write the formulas defining it within $GL(m,R)$) hence $G$ is now a type-definable subgroup of a semialgebraic group of the same dimension.\qed

\bibliographystyle{plain}

\end{document}